%% file: inoc09_wifi.tex
\documentclass[10pt]{article}
\usepackage{mathprog}
\pdfoutput=1

\usepackage{amsmath, amssymb}
\usepackage{graphics}
\usepackage{graphicx}
\usepackage{ifthen}
\usepackage{pstricks,pst-node,pst-text,pst-3d}
\usepackage{epsfig}
\usepackage{mathenv}
\usepackage{amsfonts}
\usepackage{psfrag}
\usepackage{subfigure}
\usepackage{amsopn}
\usepackage{amstext}
\usepackage{amscd}
\usepackage{latexsym}
\usepackage{multirow}
\usepackage{xspace}
\usepackage{rotating}
\usepackage{lscape}
\input epsf

\newcommand{\WPP}   {\mbox{\sc WPP}\xspace}
\newcommand{\SFAP}  {\mbox{\sc SFAP}\xspace}
\newcommand{\PSAP}  {\mbox{\sc PSAP}\xspace}
\newcommand{\CSAP}  {\mbox{\sc CSAP}\xspace}

\newcommand{\PSAPH} {\mbox{\sc PSAP-H}\xspace}
\newcommand{\PSAPL} {\mbox{\sc PSAP-L}\xspace}
\newcommand{\LinA}  {\mbox{\sc PSAP-H$_1$}\xspace}
\newcommand{\LinB}  {\mbox{\sc PSAP-H$_2$}\xspace}

\newcommand{\WFAP}  {\mbox{\sc WFAP}\xspace}

\newcommand{\WFAPH} {\mbox{\sc WFAP-H}\xspace}
\newcommand{\WFAPL} {\mbox{\sc WFAP-L}\xspace}

\newlength{\miadimensione}

\newcommand{\comment}[1]{}

\renewcommand{\leq}{\leqslant}
\renewcommand{\geq}{\geqslant}

\newcommand{\numcolumns}{\the\tab@columns}

\newcommand{\myGamma}[2]{\Gamma_{#1#2}}
\newcommand{\st}{:}

\makeatother

\newcommand{\AND}{\ \wedge\ }
\newcommand{\OR}{\ \vee\ }

\newcommand{\NP}{\mbox{\sc NP}\xspace}

\newcommand{\sss}[1]{{\scriptscriptstyle #1}}
\newcommand{\mc}[1]{\mathcal{#1}}

\newcommand{\Jij}[2]{J_{#1#2}\xspace}
\newcommand{\cJij}[2]{\tilde J_{#1#2}\xspace}

\newcommand{\CCS} {{\textrm{\tiny \sc CS}}}
\newcommand{\SF} {{\textrm{\tiny \sc SF}}}
\newcommand{\PCS} {{\textrm{\tiny \sc PS}}}

\newlength{\mytbspace}
\newcommand{\allcols}{$\times$}%

\newcommand{\tablecommands}[1]{%
\renewcommand{\allcols}{\hspace{-0.5em}$\times$}%
\fontsize{6}{7}\selectfont%
\newcommand{\mybf}{\fontsize{7}{7}\fontseries{b}\selectfont}%
\setlength{\tabcolsep}{0.3em}%
\setlength{\extrarowheight}{0.14em}%
\setlength{\mytbspace}{-0.3em}
\newcommand{\mytbspc}{\hspace{0.2cm}\ }%
\newcommand{\sol}   {\ctr{sol}}%
\newcommand{\gap}   {gap}%
\newcommand{\lpgap} {lp-gap}%
\newcommand{\hgap}  {h-gap}%
\newcommand{\cov}   {cov}%
\newcommand{\perc}  {(\%)}%
\newcommand{\Time}  {time}%
\newcommand{\Secs}  {(sec)}%
\newcommand{\of}    {(o.f.)}%
\newcommand{\groupline}{\cline{1-#1}}%
\newcommand{\lowden}{%
	\mcl{1}{l}{} \\[-0.6em]%
	\mcl{#1}{c}{\sc Low density geometric instances} \\ \groupline%
}%
\newcommand{\highden}{%
	\mcl{1}{l}{} \\[-0.6em]%
	\mcl{#1}{c}{\sc High density geometric instances} \\ \groupline%
}%
\def\quadraticlegenda{
\mcl{1}{r}{} \\
\mcl{#1}{c}{
\begin{tabular}{crl}
    gap&:    & gap between the solution value and the upper bound \\
    h-gap&:  & gap between the solution value, evaluated with the hyperbolic o.f., \\
         &   & and the best known solution for the hyperbolic o.f. \\
\end{tabular}
\hspace{3em}
\begin{tabular}{crl}
    $*$&: & gap equal to zero \\
    $-$&:  & time limit exceeded \\
\end{tabular}
} \\
\mcl{1}{r}{} \\
}%
}

\newcommand{\ctr}[1]{\hfil #1 \hfil}

\newcommand{\mcl}[3]{\multicolumn{#1}{#2}{#3}}

\newenvironment{keywords}{\vspace{0.2in} \begin{quote} \small \em {\bf Keywords\/}:}{\end{quote}}

\begin{document}

{
\centerline{\Large {\bf Modeling and Solving AP Location and Frequency Assignment}}
\centerline{\Large {\bf for Maximizing Access Efficiency in Wi-Fi Networks}}
}

\hfill

\hfill

\hfill

\centerline{{\large Sandro Bosio$^{\ast}$ \ \ Di Yuan$^{\diamond}$}}
       
\hfill

\centerline{{\small {\it ${}^{\ast}$Institut f\"ur Mathematische Optimierung, Otto-von-Guericke Universit\"at}}}
\centerline{{\small {\it 39106, Magdeburg, Germany}}}

\medskip

\centerline{{\small {\it ${}^{\diamond}$Department of Science and Technology, Link\"{o}ping University}}}
\centerline{{\small {\it SE-601 74, Norrk\"oping, Sweden}}}

\hfill 

\hfill

\noindent \hrulefill

\begin{abstract}
In this paper, we present optimization approaches for Access Point
(AP) location and frequency assignment, two major planning tasks in
deploying Wi-Fi networks.  Since APs are relatively cheap, the major
concern is network performance. We consider a performance metric,
referred to as access efficiency, that captures key aspects of how
user devices share access to the wireless medium. We propose a
two-step approach to deal with AP location and frequency assignment in
maximizing access efficiency. A novelty of our modeling approach is to
estimate, in the first step of AP location, the impact of expected
frequency availability on frequency assignment.  For each of the two
steps we derive hyperbolic formulations and their linearizations, and
we propose a promising enumerative formulation. Sample results are
reported to show the applicability of the approach.
\end{abstract}


\noindent \hrulefill

\author{
\begin{minipage}{0.5\textwidth}
	\centering
	Sandro Bosio
	\vspace{2mm} \normalsize
	\\  Institut f\"ur Mathematische Optimierung
	\\  Otto-von-Guericke Universit\"at
	\\  39106, Magdeburg, Germany
	\\  Email: bosio@mail.math.uni-magdeburg.de
\end{minipage}
\begin{minipage}{0.5\textwidth}
	\centering
	Di Yuan
	\vspace{2mm} \normalsize
	\\ Department of Science and Technology
	\\  Link\"{o}ping University
	\\ SE-601 74 Norrk\"oping, Sweden
	\\ Email: diyua@itn.liu.se
\end{minipage}
}

\section{Introduction}
\label{sec:introduction}

A Wi-Fi network consists in a set of Access Points (APs) connected to
a wired backbone. Two major planning tasks are to decide, given a set
of {\em candidate sites} (CSs) and a set of frequencies, where to
install APs and which frequency to assign to each of them.  A typical
constraint is to require coverage for a set of {\em test points} (TPs)
representing user locations, as in the classical Set Covering Problem.
As APs are relatively cheap. The major concern is network
performance.  If a TP is covered by multiple APs, a user device at the TP
will select one of them and {\em associate} to it, choosing one of the
possible transmission rates (e.g. 6, 9, 12, 18, 24, 36, 48, and 54
Mbps in IEEE 802.11g).  As rate selection depends on signal quality,
devices typically associate to the AP providing the highest signal
quality, and thus rate.

The Wi-Fi medium access mechanism (CSMA/CA) is a carrier sense
protocol: Before transmitting, a device must sense the channel as
idle.  In IEEE 802.11b and 802.11g there are up to 14 overlapping
channels, but due to inter-channel interference it is a common
practice to restrict channel assignment to the three non-overlapping
channels $1$, $6$ and~$11$.  Under this assumption, medium contention
can take place only between devices operating on the same frequency.  We
use the term {\em direct interference} to refer to a medium contention
occurring between two devices because their signals directly reach
each other.
{\em Indirect interference} refers to a signal collision at one of the associated
APs.  One example is the {\em hidden terminal} scenario: Two devices that cannot
sense each other wish to transmit to the same AP.  Carrier sensing
allows both transmissions, resulting in a collision at the AP. Note that even if
the two devices are transmitting to different APs a collision may still occur,
if the transmission from a device to
its AP reaches also the second AP, which is receiving from the other device.

To summarize, two user devices working on the same frequency
interference with each other if they are direct interferers, or if (at
least) one is associated to an AP reached by the signal of the other.
In Section \ref{sec:efficiency} we present a performance metric for
devices, referred to as {\em access efficiency}, defined as the
transmission rate scaled by the probability of successful
transmission, determined by the number of direct and indirect interferers.
For simplicity, direct interference will not be considered in this work. Its inclusion 
into models and algorithms is straightforward (see~\cite{BoYu09}).

A simplified version of this metric
has been studied and validated in~\cite{BoCaCe05}, and theoretically
and algorithmically investigated in~\cite{AmBoMa06, AmBoMaYu08}.
This paper extends these works
to consider TP-AP association, rates, and multiple frequencies.
In Sections~\ref{sec:aplocation} and~\ref{sec:frequencyassignment} we
present integer programming models for access efficiency maximization
in AP location and frequency assignment, and discuss solution
approaches.  The difference of this modeling approach with respect to others (see
e.g.~\cite{abusubaih06,Siom05,LuJRGoVa06}) is the aspect that the performance metric
originates from the Wi-Fi medium access scheme.  For a discussion on
alternative modeling and solution approaches for Wi-Fi network design,
see~\cite{BoEiGeSiYu08} and the references therein.

\section{Access Efficiency and Wi-Fi Planning}
\label{sec:efficiency}

\begin{figure}[h]
\centering
\vspace{-0.1cm}
\subfigure[\ ]{\includegraphics[height=3.8cm]{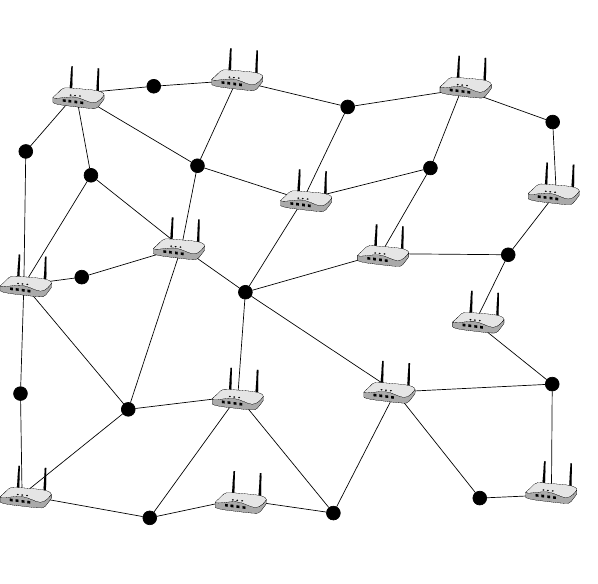}    \label{fig:covrel:a}} \hspace{1em}
\subfigure[\ ]{\includegraphics[height=3.8cm]{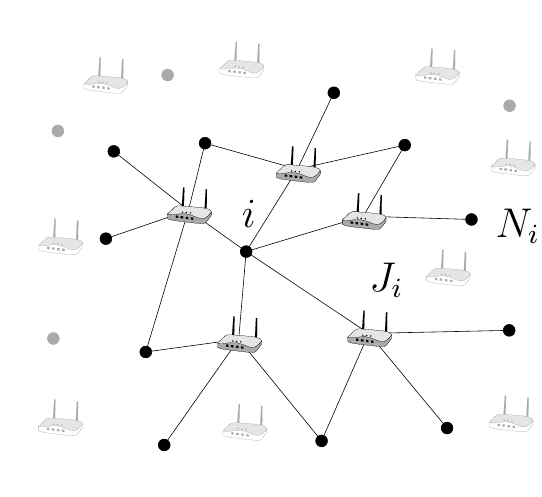}  \label{fig:covrel:b}} \hspace{1em}
\subfigure[\ ]{\includegraphics[height=3.8cm]{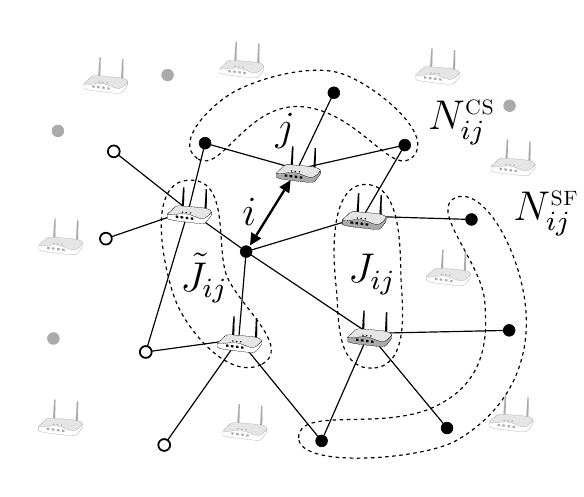} \label{fig:covrel:c}}
\vspace{-0.1cm}
\caption{
	\small \subref{fig:covrel:a} An instance,
	\subref{fig:covrel:b} covering APs and neighbors of a TP,
	\subref{fig:covrel:c} association-related sets.
}
\label{fig:covrel}
\end{figure}

The sets of TPs and CSs are denoted by $I$ and $J$ respectively, and
the set of frequencies by $F$.  We use ``AP~$j$'' to refer to an AP
installed at CS $j$, and we refer equivalently to ``TPs'' or
``users''.  The set of TPs covered by AP~$j$ and the set of APs
covering TP~$i$ are denoted respectively by $I_j \subseteq I$ and $J_i
\subseteq J$.
The set of TPs covered by a subset $S \subseteq J$ of CSs
is denoted by $I(S) = \bigcup_{j \in S}{I_j}$,
and if $I(S)=I$ the set $S$ is said to be a \emph{cover}.
The set of TPs $N_i = I(J_i)$
sharing a covering AP with TP~$i$ is referred to as the set of
\emph{neighbors} of $i$.  Figures~\ref{fig:covrel:a}
and~\ref{fig:covrel:b} illustrate the coverage and neighborhood
relations.

A design solution is a pair $(S,f)$, where $S \subseteq J$ is a cover
and $f : S \rightarrow F$ is the frequency assignment.  To model TP-AP
association, we introduce for each TP~$i$ a total order $>_i$ on
$J_i$, where $j >_i k$ means that at TP~$i$ the signal of AP~$j$ is stronger
than that of AP~$k$.  Given a cover~$S$, we denote by $a_i=a_i(S)$ the AP to which TP~$i$
associate. This is the (unique) AP $j \in S \cap J_i$ for which $j >_i
k$ for every $k \in S \cap J_i \setminus \{j\}$.
If TP~$i$ associates to AP~$j$
(shortly, if $a_i=j$) then no AP covering~$i$ with signal better than~$j$
is installed.  We denote by $\Jij{i}{j} = \{ k \in J_i : j >_i k \}$ the set of
APs covering~$i$ and compatible with the association $a_i = j$, and by
$\cJij{i}{j} = \{ k \in J_i : k >_i j \}$ those that are not compatible.  Moreover,
we define the sets of TPs $N^\CCS_{ij} = I_j \setminus \{i\}$ and $N^\SF_{ij} = I(\Jij{i}{j}) \setminus
I_j$. These notations are illustrated in Figure~\ref{fig:covrel:c}.

\vspace{-2mm}

A design solution $(S,f)$ induces a frequency assignment $f : I(S) \rightarrow F$ to the covered
TPs, with $f(i) = f(a_i)$.
Given a design solution $(S,f)$, the set
of interferers to TP $i$ can then be defined as $\Phi_i(S,f) = \{ h \in N_i : f(i)
= f(h) \AND ( a_i \in J_h \OR a_h \in J_i ) \}$.
Denoting
by $\myGamma{i}{j}$ the data rate between TP~$i$ and AP~$j$,
the access efficiency metric, under uniform traffic and equal access,
reads:

\begin{equation}
e(S,f)= \sum_{i \in I(S)} \frac{\myGamma{i}{a_i}}{1+|\Phi_i(S,f)|}.
\label{eq:eS}
\end{equation}
The {\bf Wi-Fi Planning Problem (\WPP)} consists in finding a design solution ($S$, $f$) 
maximizing $e(S,f)$.

\WPP is \NP-hard, even in the special cases $|F|=1$
(all APs working at the same frequency) and $|F|=|S|$ (complete
separation). Moreover, the frequency
assignment problem, obtained from \WPP by fixing the AP location $S$,
is also \NP-hard.
We refer to~\cite{BoYu09} for complete proofs.
In this paper we propose a
decomposition approach to deal with AP location and frequency
assignment in two steps.
A novelty of our approach is that the two tasks are not completely separated, as
in the AP location step we try to estimate the impact of frequency availability
on access efficiency.

\section{AP Location}
\label{sec:aplocation}

Consider the extreme case of planning AP location with a single
frequency (SF), i.e., $|F|=1$. The set of interfering TPs of TP~$i$, given
a cover $S$, is $\Phi_i^\SF(S) = \{ h \in N_i : a_i \in J_h \OR a_h
\in J_i \}$.  The opposite assumption is that there will be
enough frequencies to provide a frequency assignment in which two TPs interfere if and only if they are
associated to the same AP, e.g., $|F|=|S|$.  Note, however, that the
actual number of frequencies required for complete separation (CS) is
typically much smaller than $|S|$, because frequencies can be reused.  The
set of interfering TPs of TP~$i$ is then $\Phi_i^\CCS(S) = \{ h \in
N_i : a_i = a_h \}$. The access efficiency of {\bf Single Frequency AP
location} (\SFAP) and {\bf Complete Separation AP location} (\CSAP)
are respectively

\begin{equation}
e^\SF(S) = \sum_{i \in I(S)} \dfrac{\myGamma{i}{a_i}}{1+|\Phi_i^\SF(S)|} 
\mbox{~~~~~and~~~~~~}
e^\CCS(S) = \sum_{i \in I(S)} \dfrac{\myGamma{i}{a_i}}{1+|\Phi_i^\CCS(S)|}.
\label{eq:eS:SFCS}
\end{equation}

For real-life Wi-Fi deployment, the \SFAP assumption is too
conservative, while the \CSAP one is too optimistic. We therefore use
a convex combination of the two.  Consider a TP~$i \in I$ and a
cover~$S$. In \SFAP a TP $h \in N_i$ interferes with~$i$ if $a_i
\in J_h$ or $a_h \in J_i$, while in \CSAP
it interferes if $a_i = a_h$. Clearly $\Phi_i^\CCS(S)
\subseteq \Phi_i^\SF(S)$.  Users in $\Phi_i^\CCS(S)$ always interfere
with user $i$, while interference with users in $\Phi_i^\SF(S)
\setminus \Phi_i^\CCS(S)$ depends on frequency assignment.  We
consider the convex combination $\alpha |\Phi_i^\SF(S)| +
(1-\alpha) |\Phi_i^\CCS(S)|$, 
where $\alpha \in [0,1]$ represents
the likelihood that users in $\Phi_i^\SF(S) \setminus \Phi_i^\CCS(S)$
will interfere with $i$ under optimal frequency assignment.
The objective of the {\bf Partial Separation AP location Problem} (\PSAP), 
which generalizes both \SFAP and \CSAP, is to maximize

\vspace{-2mm}
\begin{equation}
	e^\PCS(S,\alpha) = 
	\sum_{i \in I(S)} 
	\dfrac{
		\myGamma{i}{ a_i}
	}{
		1+\alpha |\Phi_i^\SF(S)| + (1-\alpha) |\Phi_i^\CCS(S)|
	}
	=\sum_{i \in I(S)} 
	\dfrac{
		\myGamma{i}{ a_i}
	}{
		1 + \alpha |\Phi_i^\SF(S) \setminus \Phi_i^\CCS(S)| +
		|\Phi_i^\CCS(S)|
	}.
\label{eq:eS:PCS}
\end{equation}

\paragraph{An Integer Hyperbolic Model.}
Using the binary variables $x_j, j \in J$ for AP location,
$l_{ij}, i \in I, j \in J_i$ for TP-AP association, and $y_{ih}, i \in I, h
\in N_i, i\not=h$ for representing whether or not two TPs interfere, we obtain
a nonlinear integer model of \PSAP in which the objective functions is a
hyperbolic sum.

\vspace{-3mm}
\begin{system}{\PSAPH}{\max}{
\sum_{i \in I} \sum\limits_{j \in J_i} 
\frac{
	\myGamma{i}{j} l_{ij}
}{
	1 + \alpha \sum\limits_{h \in N_i} y_{ih} + (1-\alpha) \sum\limits_{h \in N_i} l_{hj}
}
\label{eq:PSAPH:obj}
}
\vspace{-1.5mm}
\constraint{	\sum_{j \in J_i} l_{ij} = 1
}{
		i \in I
		\label{eq:PSAPH:covering}
}%
\vspace{-1.5mm}
\constraint{	l_{ij} \leq x_j
}{
		i \in I, j \in J_i
		\label{eq:PSAPH:tocovering}
}%
\vspace{-1.5mm}
\constraint{	x_j + \sum_{k \in \Jij{i}{j}} l_{ik} \leq 1
}{
		i \in I, j \in J_i
		\label{eq:PSAPH:tonearest}
}%
\vspace{-1.5mm}
\constraint{	y_{ih} \geq \sum_{j \in J_i \cap J_h} l_{ij}
}{
		i \in I, h \in N_i
		\label{eq:PSAPH:interf1}
}%
\vspace{-1.5mm}
\constraint{	y_{ih} = y_{hi}
}{
		i \in I, h \in N_i, i<h
		\label{eq:PSAPH:symmetry}
}%
\vspace{-1.5mm}
\constraint{	x_j,l_{ij},y_{ih} \in \{0,1\}.
}{
		\nonumber }
\end{system}%
The meaning of \eqref{eq:PSAPH:covering} and \eqref{eq:PSAPH:tocovering} is obvious.
A constraint of \eqref{eq:PSAPH:tonearest} states that, if AP $j$ is installed and
covers $i$, then TP $i$ will not be associated to any AP covering $i$ with weaker
signal. By constraints \eqref{eq:PSAPH:interf1}, two TPs $i$ and $h$ interfere
($y_{ih}=1$) if at least one of them associates to an AP covering the other.
Constraints~\eqref{eq:PSAPH:symmetry} state the symmetry of the interference
relation (half of the $y$ variables can thus be removed).

A natural approach to solving \PSAPH is to
derive an integer linear reformulation. For $i \in I$ and $j \in J_i$,
we introduce a new variable $c_{ij}$ to represent the corresponding
hyperbolic term in \eqref{eq:PSAPH:obj}. The new objective function is
thus $\max \sum_{i \in I} \sum_{j \in J_i} c_{ij}$.
The value of $c_{ij}$ is defined by the nonlinear constraint
$c_{ij} + \alpha \sum_{h \in N_i} c_{ij}y_{ih} + (1-\alpha)\sum_{h \in N_i} c_{ij} l_{hj}
= \myGamma{i}{j} l_{ij}$.
We can then linearize $c_{ij}y_{ih}$ and $c_{ij}l_{hj}$ with standard techniques.
Consider for example $c_{ij}y_{ih}$.
We
introduce a new variable $z_{ihj} = c_{ij}y_{ih}$.  As $c_{ij}$ represents access
efficiency, we can derive lower and upper bounds of the value of $c_{ij}$, denoted
respectively ${\underline c}_{ij}$ and ${\bar c}_{ij}$. Then $z_{ihj} = c_{ij}y_{ih}$
can be defined by the inequalities $c_{ij} - {\bar c}_{ij} (y_{ih}-1) \leq z_{ihj}
\leq c_{ij} + {\underline c}_{ij} (y_{ih}-1)$ and ${\underline c}_{ij} y_{ih} \leq
z_{ihj} \leq {\bar c}_{ij}$.
The term $c_{ij} l_{hj}$ can be linearized in the same way.
%
%

We remark that alternative hyperbolic formulations can be obtained by
changing the definition of the $y$-variables (e.g., define
$y$-variables for $i \in I$ and $h \in \Phi_i^\SF(S) \setminus
\Phi_i^\CCS(S)$). These formulations perform similar to \PSAPH 
from a computational viewpoint.

\paragraph{An Enumerative Integer Programming Model.}

A drawback of \PSAPH is that the continuous relaxation is very weak.  A novel
integer linear model can be derived by enumerating interference scenarios for each
TP.
Recall that $\Jij{i}{j}$ and $N^\SF_{ij}$ are the sets of potential interferers
in \CSAP and \SFAP, respectively (Figure~\ref{fig:covrel:c}). We
define an {\em interference scenario} of TP $i$ as a tuple $s=(j_s, H_s, U_s)$,
where $j_s \in J_i$ is the AP to which $i$ associates, $H_s \subseteq N^\CCS_{ij_s}$
is the set of users associated to $j_s$, and $U_s \subseteq N^\SF_{ij_s}$ is the
set of users associated to some AP in $J_{ij_s}$ but not covered by $j_s$.  We
denote by $\mc S_i$ the set of all interference scenarios of TP~$i$. 
Defining new binary variables $w_{is}, i \in I, s \in \mc S_i$ to identify
interference scenario, and reusing the $x$- and $l$-variables in
\PSAPH, we obtain the following integer linear formulation.
\setlength{\extraconstraintspace}{-0.3em}
\begin{system}{\PSAPL}{\max}{
\sum_{i \in I} \sum_{s \in \mc S_i} 
\dfrac{\Gamma_{ij_s}}{1 + \alpha (|U_s| + |N^\CCS_{ij_s}|) + (1-\alpha) |H_s|}
w_{is}
\nonumber
}
\constraint{
		\eqref{eq:PSAPH:covering},\eqref{eq:PSAPH:tocovering},\eqref{eq:PSAPH:tonearest}
}{
		\nonumber
}%
\constraint{
		\sum_{s \in \mc S_i: j=j_s} w_{is} = l_{ij}
}{
		i \in I, j \in J_i \label{eq:LS:i-assignment}
}
\constraint{
		 \sum_{s \in \mc S_i: h \in U_s} w_{is} +
		 \hspace{-0.5em}
		 \sum_{j \in J_i \cap J_h} l_{ij}
 		 =
		 \hspace{-0.5em}
		 \sum_{s \in \mc S_h: i \in U_s} w_{hs} +
		 \hspace{-0.5em}
		 \sum_{j \in J_i \cap J_h} l_{hj}
}{
		i \in I, h \in N_i \label{eq:LS:ih-symmetry}
}
\constraint{
		\sum_{s \in \mc S_i: j=j_s, h \in H_s} w_{is} \leq l_{hj}
}{
		i \in I, j \in J_i, h \in N^\CCS_{ij} \label{eq:LS:ihj:leq}
		\hspace{-1em}
}
\constraint{
		\sum_{s \in \mc S_i: j=j_s, h \in H_s} w_{is} \geq l_{ij} + l_{hj} - 1
}{
		i \in I, j \in J_i, h \in N^\CCS_{ij} \hspace{2em}\label{eq:LS:ihj:geq}
		\hspace{-1em}
}
\vspace{-1.5mm}
\constraint{
		x_j, l_{ij}, w_{is} \in \{0,1\}.
}{
		\nonumber
}%
\end{system}%

The correctness of \eqref{eq:LS:i-assignment} 
is apparent.
To see the correctness of~\eqref{eq:LS:ih-symmetry}, let us consider it case by
case. For simplicity, let $W_i$, $L_i$, $W_h$ and $L_h$ represent the four summations
in~\eqref{eq:LS:ih-symmetry} (in the order they appear).  In the first case, each
of TPs~$i$ and~$h$ is associated to an AP covering also the other TP (i.e., $a_i,a_h
\in J_i \cap J_h$).  Then $L_i=L_j=1$, and $W_i=W_h=0$.  Note that also $W_i=W_h=1$
is a feasible but non-optimal assignment.  In the second case, $a_i \notin J_h$,
and $a_h \in J_i$ (or vice versa).  Then $L_i=0$ and $L_h=1$, and necessarily
$W_i=1$ and $W_h=0$.  The last case is when $a_i,a_h \notin J_i \cap J_h$, meaning
that $L_i=L_h=0$, and thus $W_i=W_h=0$.  Again, $W_i=W_h=1$ is feasible but
non-optimal.  Constraints~\eqref{eq:LS:ihj:leq} and~\eqref{eq:LS:ihj:geq} are the
linearization of the bilinear constraint $\sum_{s \in \mc S_i: j=j_s, h \in H_s}
w_{is} = l_{ij} l_{hj}$, which corresponds to the definition of interference in
\CSAP.
Note that \PSAPL can also be obtained with Dantzig-Wolfe decomposition
from an appropriate reformulation of \PSAPH.


We empirically observe that the continuous relaxation of \PSAPL yields much shaper
bounds than that of the linearization of \PSAPH.  For small or sparse instances,
one can generate the interference scenarios in advance and apply a standard solver
to \PSAPL.  For large-scale instances, the exponential number of $w$-variables can
be handled solving the LP relaxation of \PSAPL with column generation techniques.
We remark that the pricing problem decomposes by TP, and for each TP $i$ the optimal
value (i.e., highest reduced cost of $w_{is}$ among all $s \in \mc S_i$) can be found in
polynomial time. Moreover, combining column generation with branching over the
$l$-variables leads to a branch-and-price algorithm for $\PSAPH$. We omit the
details of the algorithm because of lack of space.

\section{Frequency Assignment}
\label{sec:frequencyassignment}

Given a solution $S \subseteq J$ of the AP location problem, the {\bf Wi-Fi Frequency
Assignment Problem} (\WFAP) amounts to finding a frequency assignment~$f:S \rightarrow
F$ that maximizes the efficiency $e(S,f)$.  As $S$ is given, the association between
TPs and APs is fixed.  We assume that for each AP $j \in S$ there is at least one
TP $i \in I_j$ for which $a_i = j$ (as otherwise the AP can be removed).

A fixed association corresponds to a specific scenario $s \in \mc S_i$ for each
TP $i$, which we denote by $s_i=(a_i, H_i, U_i)$, where $H_i = \{ h \in N_i \st
a_h = a_i\}$ are the neighbors of $i$ associated to $a_i$ and $U_i = \{ h \in N_i
\setminus I_{a_i} \st a_h \in J_i \}$ are the neighbors of $i$ that are not covered
by $a_i$ and are associated to some AP $a_h$ covering $i$.  Let us also denote by
$\overline U_i = U_i \cup (I_{a_i} \setminus H_i \setminus \{i\})$ the set of
neighbors of $i$ that are associated to some AP $a_h \neq a_i$ covering $i$.  Note
that $\overline U_i$ defines the set of those TPs that may or may not interfere
with $i$, depending on the frequency assignment.  The set $\Phi_i(S,f)$
of interfering TPs for a given TP~$i$ can then be more explicitly defined as
$\Phi_i(S,f) = H_i \cup \{ h \in \overline U_i \st f(i) = f(h) \}$, where the union
is disjoint by definition.
Introducing the frequency assignment
variables $x_{jf}$ ($x_{jf}=1$ if $f(j)=f$, and $0$ otherwise) and reusing 
$y$-variables, \WFAP can be formulated as follows:
\setlength{\extraconstraintspace}{-0.6em}
\begin{system}{\WFAPH}{\max}{
\sum_{i \in I} \frac{
	\myGamma{i}{a_i}
}{
	1+|H_i| + \sum\limits_{h \in \overline U_i} y_{ih}
}
\nonumber
}
\constraint{	\sum_{f \in F} x_{jf} = 1 
}{
		j \in S \label{eq:FAH2:frequency}
}
\constraint{	y_{ih} = \sum_{f \in F} x_{a_if} x_{a_hf}
}{
		i \in I, h \in \overline U_i \label{eq:FAH2:interf1}
}
\constraint{	x_{jf}, y_{ih} \in \{0,1\} .
}{
}
\end{system}%

An alternative hyperbolic formulation, based on an implicit partitioning of the
APs into $|F|$ sets, can be obtained introducing variables $v_{jk}$ ($v_{jk}=1$
if $f(j)=f(k)$, and $0$ otherwise) for all $(j,k) \in S^2$, where $S^2$ is the set
of all ordered pairs $(j,k)$ with $j,k \in S$ and $j<k$. Consider the case $|F|=3$
(a similar formulation can be derived also for $|F|=2$).
As it is easy to note that
$y_{ih}=v_{a_ia_h}$, we can write the following equivalent formulation:
\begin{system}{\WFAPH_2}{\max}{
\sum_{i \in I} \frac{
	\myGamma{i}{a_i}
}{
	1+|H_i| + \sum\limits_{h \in \overline U_i} v_{ih}
}
\nonumber
}
\constraint{	\sum_{(j,k) \in T^2} v_{jk} \geq 1
}{
		T \subseteq S \st |T|=4 \label{eq:FAH2:outgoing}
}
\constraint{	v_{jk} \geq v_{jl} + v_{kl} - 1
}{
		(j,k),(k,l) \in S^2  \label{eq:FAH2:transitivity:1}
}
\constraint{	v_{kl} \geq v_{jk} + v_{jl} - 1
}{
		(j,k),(k,l) \in S^2  \label{eq:FAH2:transitivity:2}
}
\constraint{	v_{jl} \geq v_{jk} + v_{kl} - 1
}{
		(j,k),(k,l) \in S^2  \label{eq:FAH2:transitivity:3}
}
\constraint{	v_{jk} \in \{0,1\}
}{
		(j,k) \in S^2 \nonumber,
}
\end{system}%
where $v_{ih} = v_{a_ia_h}$ if $(a_i,a_h) \in S^2$, and $v_{ih} = v_{a_ha_i}$ if $(a_h,a_i) \in S^2$.

\paragraph{An Enumerative Integer Programming Model.}

In order to define the objective function contribution for~$i$, we need to know
which TPs in $\overline U_i$ do interfere with~$i$, that is, which
APs in $C_i = \{ j \in S \st j = a_h \mbox{ for some } h \in \overline U_i\}$
operate on the same frequency of $a_i$.  By enumerating on the set $\mc A_i = \{A
\subseteq C_i\}$ of all subsets of $C_i$ we can derive the following
integer linear formulation.
\begin{system}{\WFAPL}{\max}{
\sum_{i \in I} \sum_{A \in \mc A_i} 
\dfrac{
	\myGamma{i}{a_i}
}{
	1+|H_i| + |\overline U_i \cap I(A)|
}
w_{i\sss A}
\nonumber
}
\constraint{	\eqref{eq:FAH2:outgoing},\eqref{eq:FAH2:transitivity:1},\eqref{eq:FAH2:transitivity:2}, \eqref{eq:FAH2:transitivity:3}
}{
		\nonumber
}
\constraint{	v_{a_i k} = \sum_{A \in \mc A_i \st k \in A} w_{i\sss A}
}{
		i \in I, k \in C_i     \label{eq:FALS2:interf1}
}
\constraint{	v_{jk} \in \{0,1\}, w_{ia} \in \{0,1\} ,
}{
		\nonumber
}
\end{system}%
%
It is interesting
to remark that the size of \WFAPL can be reduced by exploiting three-coloring
heuristics (assuming three frequencies) on an appropriate
\emph{AP overlap graph} of the instance, of which the node set is $S$  and the
edge set is defined as $E = \{\{j,k\} \subseteq S \st a_i \in \{j,k\} \mbox{ for
some } i \in I_j \cap I_k \}$.  If the overlap graph can
be three-colored then it is possible to find a frequency assignment $f$ for which $e(S,f) = e^\CCS(S)$.
Although three-coloring is an \NP-complete problem,
it can be decided in polynomial time for many classes of graphs.  Moreover, there
are many polynomial algorithms providing a partial three-coloring, thus allowing
size reduction of \WFAPL.

\section{Sample results and concluding remarks}
\label{sec:results}

In this section we report on some preliminary experiments carried out on random 2D instances, with both isotropic
propagation (where the coverage area of an AP is a disks)
and random anisotropic propagation (see~\cite{AmBoMaYu08} for a description of the instances).
%
In Table~\ref{table:a} we show the results obtained with four sample
anisotropic instances with 50 CSs and 100 TPs.  For each instance
problem \PSAP  is solved for five different values of the
parameter~$\alpha$, between $0.0$  (\CSAP) and $1.0$ (\SFAP). For
each resulting AP location solution, problem \WFAP is solved for both
$|F|=2$ and $|F|=3$, to derive a complete network design.
Problem \PSAP is solved to optimality with three
different formulations: two compact linearizations (\LinA and \LinB),
derived from \PSAPH with standard techniques as mentioned in
Section~\ref{sec:aplocation}, and the enumerative linearization \PSAPL.
The solution of \WFAP is not as critical as that of \PSAP,
and we do not discuss here the computational time needed.
For a comparison of different solution methods
for \WFAP the reader is referred to~\cite{BoYu09}.
Computational times reported in Table~\ref{table:a} are in seconds on
a 1281 MHz SUN UltraSPARC-IIIi with 16 GB of RAM.

\vspace{1.5em}
\begin{table}[b!ht]
\centering
\tablecommands{7}%
\begin{tabular}{r|rrr||rrr|}
\cline{2-7}
	& \mcl{1}{c}{\LinA} 
	& \mcl{1}{c}{\LinB} 
	& \mcl{1}{c||}{\PSAPL}
	& \mcl{1}{c}{\PSAP} 
	& \mcl{1}{c}{\WFAP} 
	& \mcl{1}{c|}{\WFAP}
\\
	& \mcl{1}{c}{\Time}
	& \mcl{1}{c}{\Time}
	& \mcl{1}{c||}{\Time} 
	&  
	& \mcl{1}{c}{$|F|$=2}
	& \mcl{1}{c|}{$|F|$=3}
\\
	$\alpha$~ 
	& \mcl{1}{c}{\Secs}
	& \mcl{1}{c}{\Secs} 
	& \mcl{1}{c||}{\Secs}
	& \mcl{1}{c}{\of}   
	& \mcl{1}{c}{\of}  
	& \mcl{1}{c|}{\of}
\\ \cline{2-7}
\mcl{7}{c}{} \\[-0.7em]
\mcl{7}{c}{ \sc Instance 1} \\[0.1em] \groupline \input{table.1.tex} \groupline
\mcl{7}{c}{} \\[-0.7em]
\mcl{7}{c}{ \sc Instance 2} \\[0.1em] \groupline \input{table.2.tex} \groupline
\end{tabular}
\hfil
\begin{tabular}{r|rrr||rrr|}
\cline{2-7}
	& \mcl{1}{c}{\LinA} 
	& \mcl{1}{c}{\LinB} 
	& \mcl{1}{c||}{\PSAPL}
	& \mcl{1}{c}{\PSAP} 
	& \mcl{1}{c}{\WFAP} 
	& \mcl{1}{c|}{\WFAP}
\\
	& \mcl{1}{c}{\Time}
	& \mcl{1}{c}{\Time}
	& \mcl{1}{c||}{\Time} 
	&  
	& \mcl{1}{c}{$|F|$=2}
	& \mcl{1}{c|}{$|F|$=3}
\\
	$\alpha$~ 
	& \mcl{1}{c}{\Secs}
	& \mcl{1}{c}{\Secs} 
	& \mcl{1}{c||}{\Secs}
	& \mcl{1}{c}{\of}   
	& \mcl{1}{c}{\of}  
	& \mcl{1}{c|}{\of}
\\ \cline{2-7}
\mcl{7}{c}{} \\[-0.7em]
\mcl{7}{c}{ \sc Instance 3} \\[0.1em] \groupline \input{table.4.tex} \groupline
\mcl{7}{c}{} \\[-0.7em]
\mcl{7}{c}{ \sc Instance 4} \\[0.1em] \groupline \input{table.5.tex} \groupline
\end{tabular}
\caption{
	Sample experimental results on random 2D anisotropic instances with 50 CSs and 100 TPs.
}
\label{table:a}
\end{table}

The direct linearization of \WPP cannot be solved for  these
instances within days of computation, while the direct linearizations
\LinA and \LinB of \PSAP still remain manageable, even if they show a
strong performance degeneration towards larger values of the parameter
$\alpha$, requiring up to four hours of CPU time.  In comparison, the
time needed to solve the enumerative formulation \PSAPL is by far smaller,
at most seven minutes, and shows a uniform
behavior for different values of~$\alpha$.  Note that the much better
performance showed for $\alpha=0$ and $\alpha=1$ by all the solution
methods is due to the simpler structure of \CSAP and \SFAP, which the
formulations can advantageously use.


For these instances the best performance is
obtained by choosing a small $\alpha$ in the AP location phase, i.e.,
assuming that close-to complete separation can be achieved by
frequency assignment. With three frequencies, $\alpha \leq 0.2$ gives
the highest efficiency, while with two available frequencies the
optimal value of $\alpha$ tends to increase slightly.


\paragraph{Concluding remarks.} We have presented a two-phase approach
to Wi-Fi planning, where in the AP location phase one can take into
account the estimated impact of frequency availability on an optimal
frequency assignment. For each phase we derive hyperbolic formulations
and the corresponding linearizations, and a promising enumerative
formulation approach.  The enumerative formulations allow for exact
algorithms using branch and price, of which the implementation is
underway.  More extensive experiments and result analysis will be
carried out in forthcoming work.

\paragraph{Acknowledgments.} The authors would like to thank Edoardo
Amaldi, Antonio Capone, and Federico Malucelli, Politecnico di Milano,
for intensive collaboration on WLAN design, and an anonymous referee
for helpful comments.  The work of the second author is supported by
CENIIT, Link\"oping University.

\bibliographystyle{plain}
\bibliography{article}

%
%

\end{document}

%% file: table.1.tex
0.0 & 6.27 & 14.10 & 0.11 & 742.26 & 571.25 & 711.70 \\ 
0.2 & 119.82 & 326.82 & 248.84 & 469.01 & 571.25 & 711.70 \\ 
0.4 & 632.69 & 966.48 & 330.35 & 367.91 & 576.50 & 689.06 \\ 
0.6 & 2430.33 & 2200.62 & 286.98 & 317.00 & 544.76 & 633.78 \\ 
0.8 & 6299.12 & 4440.09 & 380.36 & 286.11 & 508.66 & 581.70 \\ 
1.0 & 13580.43 & 163.99 & 55.59 & 265.75 & 484.82 & 529.47 \\ 

%% file: table.2.tex
0.0 & 4.84 & 9.05 & 0.07 & 711.13 & 599.74 & 700.38 \\ 
0.2 & 24.52 & 104.43 & 41.71 & 481.48 & 596.18 & 695.38 \\ 
0.4 & 76.83 & 199.73 & 42.05 & 379.37 & 596.18 & 695.38 \\ 
0.6 & 396.98 & 1001.74 & 42.90 & 318.82 & 589.05 & 679.08 \\ 
0.8 & 8069.79 & 14364.77 & 55.26 & 279.32 & 510.95 & 565.56 \\ 
1.0 & 1423.05 & 164.28 & 6.52 & 256.98 & 469.76 & 486.35 \\ 

%% file: table.4.tex
0.0 & 2.92 & 12.46 & 0.08 & 725.61 & 574.70 & 681.31 \\ 
0.2 & 34.07 & 135.82 & 11.50 & 476.55 & 575.41 & 695.45 \\ 
0.4 & 140.30 & 385.08 & 14.34 & 374.58 & 559.94 & 663.66 \\ 
0.6 & 466.59 & 1596.62 & 15.08 & 315.58 & 559.61 & 650.93 \\ 
0.8 & 1712.76 & 1505.59 & 14.46 & 277.06 & 536.88 & 606.36 \\ 
1.0 & 856.31 & 70.44 & 3.47 & 251.61 & 485.71 & 512.79 \\ 

%% file: table.5.tex
0.0 & 8.55 & 12.26 & 0.08 & 720.91 & 596.81 & 720.86 \\ 
0.2 & 41.88 & 201.82 & 53.03 & 499.52 & 598.68 & 712.36 \\ 
0.4 & 140.92 & 527.72 & 74.51 & 402.39 & 596.55 & 706.36 \\ 
0.6 & 925.80 & 2284.93 & 72.39 & 342.25 & 596.55 & 696.70 \\ 
0.8 & 1829.28 & 6506.13 & 82.61 & 304.99 & 518.13 & 590.03 \\ 
1.0 & 1614.43 & 363.25 & 7.05 & 284.37 & 470.91 & 506.31 \\ 